\documentclass[11pt,a4paper]{amsart}

\usepackage{aliascnt}
\usepackage{hyperref}
\usepackage{amsfonts}
\usepackage{mathrsfs}
\usepackage{amsmath}
\usepackage{amsmath,graphics}
\usepackage{amssymb}
\usepackage{indentfirst,latexsym,bm,amsmath,pstricks,amssymb,amsthm,graphicx,fancyhdr,float,color}
\usepackage[all]{xy}
\usepackage{amsthm}
\usepackage{amscd}
\usepackage[all]{hypcap}

\newtheorem{theorem}{Theorem}[section]

\newaliascnt{lemma}{theorem}

\aliascntresetthe{lemma}

\newaliascnt{conjecture}{theorem}
\newtheorem{conjecture}[conjecture]{Conjecture}
\aliascntresetthe{conjecture}

\newaliascnt{proposition}{theorem}
\newtheorem{proposition}[proposition]{Proposition}
\aliascntresetthe{proposition}

\newaliascnt{corollary}{theorem}

\aliascntresetthe{corollary}

\newaliascnt{problem}{theorem}

\aliascntresetthe{problem}

\newaliascnt{claim}{theorem}

\aliascntresetthe{claim}

\theoremstyle{definition}

\newaliascnt{definition}{theorem}

\aliascntresetthe{definition}

\newaliascnt{example}{theorem}

\aliascntresetthe{example}

\theoremstyle{remark}

\newaliascnt{remark}{theorem}
\newtheorem{remark}[remark]{Remark}
\aliascntresetthe{remark}

\newaliascnt{remarks}{theorem}

\aliascntresetthe{remarks}

\numberwithin{equation}{section}

\numberwithin{figure}{section}

\def\lra{\longrightarrow}

\def\({$($}
\def\){$)$}

\def\rank{\text{{\rm rank\,}}}

\newcommand{\rk}{{\rm rank\,}}

\newcommand{\rev}[1]{{{#1}}}

\begin{document}

\title{Strict Arakelov inequality for a family of varieties of general type}

\author{Xin Lu}

\address{School of Mathematical Sciences, Shanghai Key Laboratory of PMMP, East China Normal University,
No. 500 Dongchuan Road, Shanghai 200241, People’s Republic of China}

\email{xlv@math.ecnu.edu.cn}

\author{Jinbang Yang}

\address{Wu Wen-Tsun Key Laboratory of Mathematics, School of Mathematical Sciences, University of Science and Technology of China, Hefei, Anhui 230026, PR China}

\email{yjb@mail.ustc.edu.cn}

\author{Kang Zuo}
\address{Institut f\"ur Mathematik, Universit\"at Mainz, Mainz, Germany, 55099}
\email{zuok@uni-mainz.de}

\thanks{This work is supported by National Natural Science Foundation of China, Grant No. 12001199,
and Sponsored by Shanghai Rising-Star Program, Grant No. 20QA1403100}





\phantomsection

\maketitle

	\begin{abstract}
Let $f:\, X\to Y$ be a semistable non-isotrivial family of $n$-folds over a smooth projective curve
with discriminant locus $S \subseteq Y$  and with general fiber $F$ of general type.
We show the strict Arakelov inequality
$${\deg f_*\omega_{X/Y}^\nu \over \rk f_*\omega_{X/Y}^\nu}
< {n\nu\over 2}\cdot\deg\Omega^1_Y(\log S),$$
for all $\nu\in \mathbb N$ such that the $\nu$-th pluricanonical linear system $|\omega^\nu_F|$ is birational.
This answers a question asked by M\"oller, Viehweg and the third named author \cite{mvz06}.
\end{abstract}

\section{Introduction}
We always work over the complex number field $\mathbb C$.
Let $Y$ be a non-singular projective curve, $X$ a projective manifold,
and let $f:X\to Y$ be a proper surjective morphism with connected general fiber $F$.
Denote by $S\subseteq Y$ the discriminant divisor of $f$, i.e. $S$ is the smallest subset of points in $Y$ such that the restricted map
	$$f: X\setminus f^{-1}(S)\lra Y\setminus S $$
	is smooth. Recall that $f$ is birationally isotrivial, if
	$X\times_Y\text{Spec\,}\overline{ \mathbb C(Y)}$ is birational to
	$F\times \text{Spec\,}\overline{ \mathbb C(Y)}.$
	Putting together results due to Parshin-Arakelov, Migliorini, Zhang, Kovacs, Bedulev-Viehweg, Oguiso-Viehweg, Viehweg-Zuo, etc.
	(see \cite{vz02} and the references given there), one has
\begin{theorem}\label{thm-1}
 Let $f:X\to Y$ be a \rev{non-birationally isotrivial} family of $n$-folds, with general fiber $F$. Assume either
	\begin{itemize}
	\item $\kappa(F)=\dim (F)$, or
	\item $F$ has a minimal model $F'$ with $\omega_{F'}$ semi-ample.
	\end{itemize}
	Then $(Y,S)$ is logarithmic hyperbolic, i.e., $\deg\Omega^1_Y(\log S)>0$.
\end{theorem}

Let us mention the following theorem by Bogomolov-B\"ohning-Graf von Bothmer \cite{bbg} characterizing birationally isotrivial families. Here we should also remark that
in their paper the field $\mathbb{C}$ can be replaced by an arbitrary algebraically closed field with infinite transcendence degree over the prime field.

\begin{theorem}[Bogomolov-B\"ohning-Graf von Bothmer]\label{thm-bbg}
Let $f: V\to U$ be a family of algebraic varieties over $\mathbb{C}$ such that all fibers are birational to each other and $U$ is integral. Then $f$ is birationally isotrivial.
\end{theorem}

Let $ M_h$ denote the coarse moduli space of polarized manifolds with semi-ample canonical line bundle and with fixed Hilbert polynomial $h. $
\autoref{thm-1} is equivalent to saying that the moduli stack of $M_h$ is algebraic hyperbolic, which can be followed from the existence of a big subsheaf $\mathcal A\hookrightarrow S^\ell\Omega^1_Y(\log S)$.
As a complex analytic version, Viehweg-Zuo have constructed a complex
Finsler metric $h_f$ with strictly negative curvature and consequently, the various complex hyperbolicities;
for example, the Brody \cite{vz03}, Kobayashi \cite{ty15} and big Picard hyperbolicities \cite{dlsz19} have been proven for the moduli stack.
We remark also that when $\kappa(F)=\dim (F)$, i.e., the fiber is of general type, \autoref{thm-1} is proved in \cite{vz01}.
In the following, we briefly explain the construction of the
big \rev{subsheaf} in $S^\ell \Omega^1_Y(\log S)$ and the Finsler metric $h_f$.
It relies on three main steps:
\begin{enumerate}
\item
Kawamata-Viehweg's positivity theorem: $A:=\det f_*\omega^\nu_{X/Y}$ is big.
\item
The deformation Higgs bundle introduced by Viehweg-Zuo \cite{vz01}:
by taking a suitable $m$-th power of the \rev{self-fiber product}
$$f^{(m)}: X^{(m)}\to Y$$
such that $A^m\subset f_*^{(m)}\omega^\nu_{X^{(m)}/Y}$
and running the maximal non-zero iteration of the Kodaira-Spencer map on the deformation Higgs bundle
$$(F,\tau)\otimes A=(\bigoplus_{p+q=mn} F^{p,q},\bigoplus_{p+q=mn}\tau^{p,q})\otimes A$$ attached to the family $f^{(m)}$ twisted with $A$, one obtains the Griffiths-Yukawa coupling:
$$ S^{\ell}T_{Y}(-\log S): A\xrightarrow{\tau^\ell} F^{mn-\ell,\ell}\otimes A,$$
which induces
an inclusion of sheaves
$$ \mathcal A:=A\otimes P\hookrightarrow S^\ell\Omega^1_Y(\log S).$$
Here $P$ is the dual of $\text{Im}(\tau^\ell)$, which is contained in the kernel of the Kodaira-Spencer map on the next graded piece
$$\tau^{mn-\ell-1, \ell+1}: F^{mn-\ell ,\ell}\otimes A\to (F^{mn-\ell-1,\ell+1}\otimes A)\otimes\Omega^1_Y(\log S).$$
We will show $P$ is semi-positive in the next step by applying Hodge theory. Hence, $\mathcal A$ is big.
\item
A comparison map: by taking the $\nu$-th cyclic cover $Z\to X^{(m)}$ defined by a section in the linear system $|\omega^\nu_{X^{(m)}/Y}\otimes f^{{(m)}*}A^{-\nu}|$ and taking the graded Higgs bundle $(E,\theta)$ as the grading of the \rev{quasi-canonical} extension of the variation of the Hodge structure on the middle cohomology attached to the induced family
$$g: Z\to X^{(m)}\xrightarrow{f}Y,$$
one constructs a comparison map of Higgs bundles
$$\rho: (F,\tau)\otimes A\to (E,\theta).$$
The semi-negativity of $\text{ker}(\theta)$ (\cite{zuo00}) shows that
$\text{ker}(\tau)$ is semi-negative. Hence, the sheaf $P$ appearing
in Griffiths-Yukawa coupling as the dual of $\text{Im}(\tau^\ell)$ contained in $\text{ker}(\tau)$ is semi-positive.
\end{enumerate}

Given a complex manifold $M$ with a Hermitian metric such that the holomorphic sectional curvature is bounded above by a negative $(1,1)$-form, Yau's Schwarz lemma says that for any holomorphic map $\gamma$ from a hyperbolic curve $ Y$ into $M$, the pull-back  metric to $Y$ is bounded from above by the hyperbolic metric on $Y$.
In our situation, let $(Y,S)\to (M, \partial M)$ be a logarithmic hyperbolic curve in the moduli stack, which carries a complex Finsler metric with the holomorphic sectional curvature bounded from above by a negative $(1,1)$-form. Then the global form of
Yau's Schwarz lemma can be expressed as an upper bound of the degree of the big subsheaf
$$ \mathcal A=A\otimes P\hookrightarrow S^\ell\Omega^1_Y(\log S)$$
in terms of $\deg\Omega^1_Y(\log S).$
As $P$ is a nef invertible sheaf, we obtain also an upper bound
 $$\deg f_*\omega^\nu_{X/Y}=\deg A|_Y\leq \ell\cdot\deg\Omega^1_Y(\log S), $$
 the so-called Arakelov inequality for the direct image of relative pluri-dualizing sheaf.

 We are more interested in this type of Arakelov inequality with an explicit and optimal upper bound.
 It is well-known for $\nu=1$ since a long time.
Faltings and Deligne have \rev{proved} 1th Arakelov inequality for families of abelian varieties.
For a semistable family $f:X\to Y$ of abelian $g$-folds, Faltings \cite{fal83} (with an improvement by Deligne \cite{del87}) showed that
$$\deg f_*\omega_{X/Y}\leq {g\over 2}\deg\Omega^1_Y(\log S).$$
As $f_*\omega^\nu_{X/Y}=(f_*\omega_{X/Y})^\nu$ for a semistable family of abelian varieties one
obtains immediately the $\nu$-th Arakelov inequality by taking $\nu$-th power of 1th Arakelov inequality
$$ \deg f_*\omega_{X/Y}^\nu=\deg (f_*\omega_{X/Y})^\nu
\leq {\nu g\over 2}\deg\Omega^1_Y(\log S).$$
  Peters \cite{per00}, Jost-Zuo \cite{jz02},
 Viehweg-Zuo \cite{vz03-2}, Green-Griffiths-Kerr \cite{ggk09}, and a very recent work by Biquard-Collier-Garcia-Prada-Toledo \cite{bcgt}
have studied Arakelov inequality for systems of Hodge bundles over curves.

Tan \cite{tan95} and Liu \cite{liu96} have shown the strict Arakelov inequality for semistable families of curves of genus $g\geq 2$, namely
$$\deg f_*\omega_{X/Y}< {g\over 2}\deg\Omega^1_Y(\log S).$$
The proof of Tan relies on the strict logarithmic Miyaoka-Yau inequality on fibred algebraic surfaces \cite{miy84} and Xiao's slope inequality \cite{xia87}, while Liu uses the explicit and optimal \rev{upper bound} of the holomorphic sectional curvature of the Weil-Peterson metric on the moduli space of curves.

In \cite{vz04} it is shown that the $\nu$-th Arakelov equality holds for
a semistable family of abelian varieties if and only if the family is a universal family over a Shimura curve of Mumford-Tate type.

Viehweg-Zuo and M\"oller-Viehweg-Zuo started to study
the Arakelov inequality for the direct image of pluri relative dualizing sheaf of a family of arbitrary fibers.
For a vector bundle $W$ over a smooth projective curve $Y$ the slope of $W$ is defined as
$$\mu(W)={\deg W\over {\rk} W}.$$
In \cite{mvz06,vz06} one finds
\begin{theorem}\label{thm-2}
 Assume that $f:X\to Y$ is a semistable family of $n$-folds over a smooth projective curve $Y$ with the discriminant locus $S$. If $Y=\mathbb P^1$ assume in addition that $\#S\geq 2$. Then for all $\nu\geq 1$ with $f_*\omega^\nu_{X/Y}\not=0$ one has
 $$\mu( f_*\omega_{X/Y}^\nu)
\leq {n\nu\over 2}\cdot\deg\Omega^1_Y(\log S).$$
More generally, for any non-zero subbundle $W\subset f_*\omega^\nu_{X/Y}$ it holds that
\begin{equation}\label{eqn-1-11}
 \mu(W)\leq {n\nu\over 2}\cdot\deg\Omega^1_Y(\log S).
\end{equation}
\end{theorem}

We remark that the above theorem recovers the Arakelov inequality for the case $\nu=1$ due to Faltings and Deligne.
Moreover, motivated by the Colemen-Oort conjecture (cf. \cite{mo13})
that the Torelli locus of curves of genus $g\gg 0$ can not contain generically any Shimura subvariety of positive dimension
and the characterization of Shimura families using the Arakelov equality by Viehweg-Zuo \cite{vz04}, we pose
\begin{conjecture}
 Let $f:X\to Y$ be a semistable non-birationally isotrivial family of $n$-folds over a smooth projective curve $Y$ with the discriminant locus $S$. Assume that the Arakelov equalities hold, i.e.,
 $$\mu( f_*\omega_{X/Y}^\nu)
= {n\nu\over 2}\cdot\deg\Omega^1_Y(\log S),\quad \forall\, \nu\in \mathbb N,\,\text{with}\, f_*\omega^\nu_{X/Y}\not=0.$$
Then the general fiber is of Kodaira dimension zero and the family is Shimura in the following sense: $Y$ parameterizes a compactified universal family of abelian varieties $f': X'\to Y$ with the given Mumford-Tate group, and the variation of Hodge structure on the middle cohomology of $f: X\to Y$ is a direct factor of a tensor product of the weight-1 variation of Hodge structure attached to the universal family of abelian varieties $f': X'\to Y $.
\end{conjecture}

 Viehweg-Zuo and M\"oller-Viehweg-Zuo \cite{mvz06,vz06}
 have generalized the strict Arakelov inequality due to Tan and Liu for the direct image of
 the relative dualizing sheaf on a semistable family of manifolds of general type:
\begin{theorem}[\cite{vz06}]\label{thm-3}
 Let $f: X\to Y$ be a semistable and {non-birationally isotrivial} family of $n$-folds, and
let $ W\subseteq f_*\omega_{X/Y}$ be a subbundle. Assume that either
\begin{itemize}
\item[\bf a.~] $f^*W\to \omega_{X/Y}$ defines a birational $Y$-morphism $\eta: X\to \mathbb P_{Y}(A),$ or
\item[\bf b.~] $n=1$ and $\rk W\geq 2$.
\end{itemize}
Then
$$\mu(W)< {n\over 2}\deg\Omega^1_Y(\log S).$$
\end{theorem}

The assumption {\bf a.} in \autoref{thm-3} seems to be a little bit too strong, as the pluri-canonical linear system $|\omega^\nu|$ on a variety of general type defines a non-constant, or birational map only after taking a higher power.
In this note we show the strict Arakelov inequality for a semistable family of $n$-folds of general type in the following form,
which answers a question asked by M\"oller, Viehweg and the third named author; see the discussion after \cite[Theorem\,0.3]{mvz06}.
\begin{theorem}\label{thm-4}
 Let $f:X\to Y$ be a non-isotrivial semistable family of $n$-folds, and let
$ W\subseteq f_*\omega^\nu_{X/Y}$ be a subbundle. Assume that $f^*W\to \omega^\nu_{X/Y}$ defines a birational $Y$-morphism $\eta: X\to \mathbb P_{Y}(W),$
then
$$\mu (W)
< {n\nu\over 2}\cdot\deg\Omega^1_Y(\log S).$$
\end{theorem}

The above theorem leads some direct consequences for families of manifolds of general type of small dimensions.\\[.2cm]
{\bf 1.} For a semistable family $f:\,X\to Y$ of minimal surfaces of general type, which is {non-isotrivial} (note that both notions birationally isotrivial and isotrivial are equivalent for a minimal surface)
then it is known that the $\nu$-th pluricanonical linear system on a general fiber defines a birational map for $\nu\geq 5$ (cf. \cite{rei88}).
Hence we obtain the strictly Arakelov inequality
$${\deg f_*\omega_{X/Y}^\nu \over {\rk} f_*\omega_{X/Y}^\nu}
< \nu\cdot\deg\Omega^1_Y(\log S),\qquad \forall~\nu\geq 5.$$\\[.2cm]
{\bf 2.} M. Chen and \rev{J. Chen} \cite{cc07} have shown for any smooth three-fold $F$ of general type there exists a number $m(3)$
between 27 and 57 and depending on certain classification classes on $F$
such that
$\omega_F^\nu$ defines a birational map for any $v\geq m(3).$ So we obtain
the strictly Arakelov inequality of a semistable family of three-folds of general type
$${\deg f_*\omega_{X/Y}^\nu \over {\rk} f_*\omega_{X/Y}^\nu}
< {3\nu\over 2}\cdot\deg\Omega^1_Y(\log S), \qquad \forall~\nu \geq 57.\vspace{1mm}$$

We would also like to point out that the Arakelov bound on the slope of subbundles $W\subset f_*\omega^\nu_{X/Y}$
is asymptotically optimal in the following sense:
there exist semistable families of $n$-folds and subbundles $W_{\nu}\subseteq f_*\omega^\nu_{X/Y}$ such that
$f^*W_{\nu}\to \omega^\nu_{X/Y}$ defines a birational $Y$-morphism $\eta: X\to \mathbb P_{Y}(W_{\nu})$, and that
\begin{equation}\label{eqn-0-1}
 \lim_{\nu \to \infty} \frac{\mu(W_{\nu})}{{n\nu\over 2}\cdot\deg\Omega^1_Y(\log S)}=1.
\end{equation}
Indeed, let $f:\,X \to Y$ be the universal family over a Teichm\"uller curve \cite{mol06}.
Then there is a line bundle $L\subseteq f_*\omega_{X/Y}$ such that
$$\mu(L)=\deg(L)=\frac12\deg\Omega^1_Y(\log S).$$
Clearly,
$$W_{\nu}:=L^{\nu-1}\otimes f_*\omega_{X/Y} \subseteq f_*\omega_{X/Y}^{\nu},\qquad \forall\,\nu\geq 2.$$
Note that $f^*W_{\nu}\to \omega^\nu_{X/Y}$ defines a birational $Y$-morphism $\eta: X\to \mathbb P_{Y}(W_{\nu})$ if $f$ is non-hyperelliptic,
because the morphism defined by $f^*W_{\nu}\to \omega^\nu_{X/Y}$ is the same as the one defined $f^*f_*\omega_{X/Y}\to \omega_{X/Y}$,
which is birational when $f$ is non-hyperelliptic.
Moreover, by direct computation, one has $\mu(W_{\nu})=(\nu-1)\mu(L)+\mu(f_*\omega_{X/Y})$.
Thus \eqref{eqn-0-1} holds.
Taking the \rev{self-fiber product}, one can get semistable families of $n$-folds with subbundles $W_{\nu}\subseteq f_*\omega_{X/Y}^\nu$ satisfying \eqref{eqn-0-1}.

The statement {\bf b.} in \autoref{thm-3} is a strong result.
The proof relies on the theory on Teichm\"uller theory due to M\"oller.
At the moment we do not know how to prove this type of strict inequality for semistable families of higher dimensional varieties.
We leave this as a conjecture.
\begin{conjecture}
 Let $f:X\to Y$ be a semistable family of $n$-folds of general type, then for any subbundle $W\subseteq f_*\omega_{X/Y}^\nu$
 of ${\rk} W\geq 2$ the following inequality holds
$$\mu(W)< {n\nu\over 2}\cdot\deg\Omega^1_Y(\log S).$$
\end{conjecture}

\vspace{2mm}
The structure of the note is organized as follows.
 In \autoref{section-2} for readers' convenience we recall the notion of the deformation Higgs bundle attached to a family, in particular, over a \rev{1-dimensional} base and we explain the comparison between the deformation Higgs bundle and the graded Higgs bundle attached to a cyclic cover
 of the original family.

 In \autoref{section-4} we first sketch the proof of \autoref{thm-2}.
 Given a subbundle $W\subset f_*\omega^\nu_{X/Y}$, by taking
 $\det W$ and a \rev{self-fiber product} of $f:X\to Y$ of a suitable power we may reduce to an invertible \rev{subsheaf} $L\subseteq f_*\omega^\nu_{X/Y}$. Furthermore, by taking a base change of $Y$
we may raise the $\nu$-power of $L$ and obtain a section
$$s: \mathcal O_X\to \omega^\nu_{X/Y}\otimes f^*L^{-\nu}.$$
 Via the new family induced by the $\nu$-cyclic cover $h: W\xrightarrow{\tau} X\xrightarrow{f} Y$ by taking the $\nu$-th roots out of the section $s,$ we construct a comparison map between the deformation Higgs bundle twisted with $L$
and the logarithmic graded Higgs bundle of the variation of Hodge structure of the middle cohomology of the new family $h$
$$\rho: (F,\tau)\otimes L\to (E,\theta).$$
By applying Simpson's Higgs semistability on the Higgs subbundle the deformation Higgs bundle in $(F,\tau)\otimes L$ generated by $L=L\otimes F^{n,0}$ and $\tau$ via the comparison map we complete the proof of \autoref{thm-2}.

\autoref{thm-4} will be proved in \autoref{sec-4-2}.
Simpson's theorem \cite{sim92} on the formality of the category of semistable vector bundles on a smooth projective curve of degree zero plays a crucial role in the proof.
Given a subbundle $W\subseteq f_*\omega^\nu_{X/Y}$ and assuming that
the \rev{sub-linear system} $f^*W\to\omega^\nu_{X/Y}$ defines a birational map
$$\eta: X\to\mathbb P_{Y}(W).$$
Then the family $f':X':=\eta(X)\to Y$ induced by the map $\eta$ is \rev{non-birationally isotrivial}. The crucial observation in the proof of \autoref{thm-4} is that
if the Arakelov inequality for $W$ becomes an equality then the family
$ f': X'\to Y$ must be birationally isotrivial.
More precisely, we consider the $d$-multiplication map
$$ S^d W\subseteq S^d f_*\omega^\nu_{X/Y}\xrightarrow {m_d} f_*\omega^{\nu d}_{X/Y}, \quad \text{for} \quad d\in \mathbb N.$$
It is well-known by the definition of a map induced by a linear system the kernel $K_d\subseteq S^d(W)$ of the multiplication map $m_d$ restricted to $S^d(W)$
is just the sheaf of \rev{homogeneous} polynomials of degree $d$ in the \rev{homogeneous} ideal defining the fibers of $f': X'\to Y.$
If $\mu(W)$ achieves the maximal value
$$\mu(W)= {n\nu\over 2}\cdot\deg\Omega^1_Y(\log S),$$
then by applying \autoref{thm-2} to \rev{subsheaves} of $W$ we know that
$W$, and hence all powers $S^d(W)$ are semistable. A simple semistability argument on the multiplication maps shows that all
kernels $K_d\subseteq S^d(W)$ remain semistable and have the same slope as $S^d(W)$'s. So,
after a base change and twisting $W$ with a line bundle, we may assume $ S^d(W) $ is semistable of degree 0 and $K_d\subset S^d(W)$ is of degree $0$ for all $d\in \mathbb N$.
By solving the (approximate) Yang-Mills equation on the semistable vector bundle $W$ of degree zero on the smooth projective curve $Y$,
one obtains integrable connections on $W$.
Thanks to Simpson's theorem (cf. \cite{sim92}) we find a canonical integrable connection $(W,\nabla^\text{can})$ in the sense that
all $K_d\subseteq S^d(W)$ are preserved by the induced connection $(S^d(W), S^d(\nabla^\text{can})).$
This means also that for each point $y\in Y$ we may find an analytic open neighborhood $y\in U\subset Y$ and find a flat base $\mathbb W$ of $(W,\nabla^\text{can})(U)$ ,
such that the flat space $\mathbb K_d(U)$ of $(K_d, S^d(\nabla^\text{can}))(U)$ is a \rev{subspace}
of $S^d(\mathbb W)$. This implies that under the flat base $\mathbb W$ of $(W_U, \nabla_U)$ and for any $d\in \mathbb N$ the coefficients of
all \rev{homogeneous} polynomials of degree $d$ in the \rev{homogeneous} ideal defining the fibers of $f': X'_U\to U$ are constant up to a scalar multiplication. This shows that
the family $f': X'_{U}\to U$ is constant.
But, it leads to a contradiction to $\eta$ is a birational embedding and our original family is \rev{non-birationally isotrivial} (cf. \autoref{thm-bbg}).

\section{Deformation Higgs Bundle and comparison
with variation of Hodge structures}\label{section-2}
\subsection{ Graded Higgs bundle arising
from variation of Hodge structures}
Throughout this section, we will assume that $U$ is a quasi-projective manifold and compactified by a projective manifold
$\overline{Y}$ with $\overline{S}=\overline{Y}\setminus U$ \rev{being} a simple normal crossing divisor,
and that there is smooth family $f: V\to U$ of $n$-folds.
Though in this note we only consider a family over a 1-dimensional base, for reader's convenience we recall some basic facts about Hodge theory attached to family of $n$-folds over a base of arbitrary dimension.\\
Leaving out a codimension two subset of $\overline{Y}$ we find a good partial compactification $f:X\to Y$ in the following sense
\begin{itemize}
\item $X$ and $Y$ are quasi-projective manifolds, $f$ is flat, $U\subseteq Y$ and
$\mathrm{codim}(\overline{Y}\setminus Y)\geq 2.$
\item $S=Y\setminus U$ is smooth and $\Delta=f^*S$ is a relative simple normal crossing divisor over $S$ (i.e. whose components, and all their intersections are smooth over components of $S$).
\end{itemize}
Following Griffiths and Simpson, one constructs the most natural graded Higgs bundle (or, system of Hodge bundles called by Simpson) related to the geometry and topology on the family.
Taking the wedge product, one sees that the tautological sequence
\begin{equation}\label{eqn-2-1}
 0\to f^*\Omega^1_Y(\log S)\to\Omega^1_X(\log\Delta)\to\Omega^1_{X/Y}(\log\Delta)\to 0
\end{equation}
induces the short exact sequence of logarithmic forms of higher degrees
\begin{equation}\label{eqn-1-2}
 0\to f^*(\Omega^1_Y(\log S))\otimes\Omega_{X/Y}^{p-1}(\log \Delta)\to\text{gr\,}\Omega^p_X(\log \Delta)\rightarrow\Omega^p_{X/Y}(\log \Delta)\to 0,
\end{equation}
where
\[\text{gr\,}\Omega^p_{X}(\log\Delta)=\Omega^p_X(\log\Delta)/f^*\Omega^2_Y(\log S)\otimes\Omega^{p-2}_X(\log\Delta).\]
The direct sum of the direct image sheaves
$$E^{p,q}=R^qf_*\Omega^p_{X/Y}(\log\Delta),\quad p+q=k$$ endowed with the connecting maps in \eqref{eqn-1-2}
\[\theta^{p,q}: R^qf_*\Omega^p_{X/Y}(\log\Delta)
\xrightarrow{\partial}\Omega^1_Y(\log S)\otimes R^{q+1}f_*\Omega^{p-1}_{X/Y}(\log\Delta)\]
forms a so-called system of Hodge bundles of weight-$k$ by Simpson.
\[(E,\theta)
=(\bigoplus_{p+q=k}E^{p,q},\bigoplus_{p+q=k}\theta^{p,q}).\]
Take the dual of\eqref{eqn-2-1}, one has an exact sequence
$$0\to T_{X/Y}(-\log\Delta)\to T_X(-\log\Delta)\to f^*T_Y(-\log S)\to 0.$$
The connecting map of the direct image defines the logarithmic Kodaira-Spencer map
\[\tau: T_Y(-\log S)\to R^1f_*T_{X/Y}(-\log\Delta).\]
The Higgs field $\theta^{p,q}$ can be also defined as the cup-product with $\tau$
\begin{equation*}
\xymatrix@R=15mm{
T_Y(-\log S)\otimes R^qf_*\Omega^p_{X/Y}(\log\Delta)
\ar[r]^-{\theta^{p,q}} \ar[d]^-{\tau\circ \mathrm{id}}
&
R^{q+1}f_*\Omega^{p-1}_{X/Y}(\log\Delta)
\\
R^1f_*T_{X/Y}(-\log\Delta)\otimes R^qf_*\Omega^p_{X/Y}(\log\Delta)
\ar[ur]^-{\cup}
&\\
}
\end{equation*}

\begin{proposition}\label{prop-1-1}
 The Higgs bundle $(E,\theta)$ is the grading of Deligne's quasi-canonical extension of the variation of the polarized Hodge structure on $k$-th Betti cohomology $R^kf_*\mathbb Z_{V}$ of the smooth family $f:V\to U$.
\end{proposition}
This result is well-known and due to Griffiths \cite{gri70}.
Katz-Oda \cite{ko68} have an algebraic approach, which works also over any characteristic satisfying $E_1$-degeneration of Hodge to de Rham spectral sequence.
The Higgs field $\theta: E\to E\otimes \Omega^1_Y(\log S)$
induces a natural map
$$\vartheta: T_Y(-\log S)\to\mathcal End(E),$$
which coincides with the derivative of the period mapping attached to the variation of Hodge structure on $\mathbb R^kf_*\mathbb Z_V$ over $U.$
By Griffiths' curvature formula
it known the lass Hodge bundle $E^{0,k}$ and, in a slightly general form, the kernel of the Higgs field is semi-negative \cite{zuo00}.
 Assume that the period mapping is locally injective
 (equivalently, $\theta$ is injective over $U$) for example, families of hypersurfaces in projective space of high degrees.
 Then by Griffiths-Schmid's theorem on the curvature of the Hodge metric along the horizontal direction
 in the period domain we know that the holomorphic sectional curvature of the pulled back Hodge metric on the base $U=Y\setminus S$ is bounded above by a negative $(1,1)$-form.
However, we notice that the Torelli injectivity could fail for general varieties.
For example, for surfaces of general type with small Chern classes the period mapping can be constant.

In the joint work \cite{vz02} Viehweg and the third named author have started looking for a replacement of variation of Hodge structure in the case when the Torelli injectivity fails, the so-called deformation Higgs bundle. We are going to
briefly discuss the construction in the next section.
\subsection{Deformation Higgs bundle attached to a family}
Given a log smooth family $f: (X,\Delta)\to (Y, S),$
we start with the classical logarithmic Kodaira-Spencer map
$$ T_Y(-\log S)\xrightarrow{\tau^{n,0}}R^1f_*T_{X/Y}(-\log\Delta).$$
The Kodaira-Spencer map measures the variation of complex structure.
The Kodaira-Spencer map $\tau^{n,0}|_U$ is zero if and only if the smooth family $f: V\to U$ is isotrivial, cf. \cite{kod}.

In a similar way as we have done for the Higgs field on a system of Hodge bundles Viehweg-Zuo \cite{vz02} introduced the extended Kodaira-Spencer map $\tau^{p,q}$ as follows:

Put $\mathcal L:=\Omega^n_{X/Y}(\log \Delta)$ and
 consider the tautological exact sequence of logarithmic forms of higher degree twisted by $\mathcal L^{-1}$
 \begin{equation}\label{eqn}
 0\to f^*(\Omega^1_Y)\otimes\Omega^{p-1}(\log S)\otimes \mathcal L^{-1}\to\text{gr\,}\Omega^p_X(\log S)\otimes \mathcal L^{-1}\rightarrow\Omega^p_{X/Y}(\log S)\otimes \mathcal L^{-1}\to 0.
\end{equation}
For $p+q=n$, we take
$$F^{p,q}:=R^qf_*(\wedge^qT_{X/Y}(-\log \Delta))/\text{torsion} = R^qf_*(\Omega^p_{X/Y}(\log\Delta)\otimes\mathcal L^{-1})/\text{torsion},$$
and define $\tau^{p,q}$ as the connecting map at the place $q\to q+1$:
$$ \tau^{p,q}: F^{p,q} \to
 F^{p-1,q+1}\otimes\Omega^1_Y(\log S).$$
Putting all individual sheaves
$F^{p,q}$
together and endowed with the maps $\tau^{p,q}$ we obtain the so-called {\sl Deformation Higgs bundle (sheaf)} attached to $f:X\to Y$:
$$(F,\tau):=\Big(\bigoplus_{p+q=n}F^{p,q},\bigoplus_{p+q=n}\tau^{p,q}\Big).$$
We remark that the extended Kodaira-Spencer map can also be represented as the cup product in a standard way
\begin{equation*}
\xymatrix@R=2cm{
T_Y(-\log S)\otimes R^qf_*T^q_{X/Y}(-\log\Delta)\ar[d]^-{\tau^{n,0}\otimes\mathrm{Id}}\ar[r]^-{\tau^{p,q}}& R^{q+1}f_*T^{q+1}_{X/Y}(-\log\Delta)\\
R^1f_*T_{X/Y}(-\log\Delta)\otimes R^qf_*T^q_{X/Y}(-\log\Delta).\ar[ur]^-{\cup}
}
\end{equation*}
The extended Kodaira-Spencer map $\tau$
 satisfies the integrability condition $\tau\wedge\tau=0$.
Indeed, for a 1-dimensional base $Y$ considered in this note
 the integrability holds trivially true.
 In general, for a higher dimensional base,
 using Dolbeault representative for $H^q(X_y, T^q_{X_y})$ the cup product $\cup$ in the above diagram
 is nothing but the usual wedge product of bundle-valued differential forms,
and the integrability just follows from the commutativity of the wedge product of differential forms of even degrees.

\subsection{Kawamata-Viehweg's positivity of direct image sheaves}
Let $f:X\to Y$ be a family of $n$-folds over a 1-dimensional base $Y$ with semistable singular fibers $\Delta$ over $S$ and with the smooth part of the family
$$f: V=X\setminus\Delta\to Y\setminus S=:U.$$
As the family is semistable $\Omega^n_{X/Y}(\log\Delta)=\omega_{X/Y},$ we recall following positivity of $f_*\omega_{X/Y}^\nu$
\begin{theorem}[Kawamata and Viehweg, cf. \cite{kaw85,vie83}]
Assume that $f: X\to Y$ has the maximal variation, and $\omega _{V/U}$ is semi-ample.
Then $f_*\omega_{X/Y}^\nu $ is weakly positive for all $\nu>1$ with $f_*\omega_{X/Y}^\nu\not=0.$
\end{theorem}

\subsection{ Comparison between deformation Higgs bundle and system of Hodge bundles}
The comparison map relies on the certain type of cyclic covers on the family $f:X\to Y$. The motivation of constructing cyclic covers
goes back to the work by Esnault-Viehweg \cite{es92}. They gave
 a more Hodge theoretical approach to the Kodaira-Akizuki-Nakano vanishing theorem.

\subsubsection{Comparison map in the absolute case and Kodaira type vanishing theorem}
Let $X$ be a projective manifold, $\mathcal L$ an ample line bundle on $X$ and $s\in H^0(X,\mathcal L^\nu)$ with the simple normal crossing zero divisor $D:=(s)_0\subseteq X.$ One takes the $\nu$-th cyclic cover
$$\gamma: Z= X(\sqrt[\nu]{s}) \to X$$
with
$$\gamma_*\Omega^p_Z(\log \gamma^*D)
=\bigoplus_{i=0}^{\nu-1}\Omega^p_X(\log D)\otimes \mathcal L^{-i}.$$
Deligne has shown
$$ H^k(Z\setminus\gamma^*D,\mathbb C)=\bigoplus_{p+q=k} H^q(Z,\Omega^p_Z(\log \gamma^*D)).$$
Assume $D$ is ample, then $X\setminus D$ is affine, the same holds true for $Z\setminus \gamma^*D$
and hence $$H^k(Z\setminus \gamma^*D, \mathbb C)=0,\quad \forall\, k>\dim X=n.$$\\
By the Hodge decomposition
$$0=H^q(Z,\Omega^p_Z(\log\gamma^* D)=
\bigoplus_{i=0}^{\nu-1}H^q(X, \Omega^p_X(\log D)\otimes\mathcal L^{-i}).$$
for any $p+q>n.$ In particular,
$$ H^q(X,\Omega^p_X(\log D)\otimes\mathcal L^{-1})=0,\forall\, p+q>n.$$
Using the residue map as well as the Serre duality
one also shows the
Kodaira-Akizuki-Nakano vanishing theorem by induction on $\dim X$:
\begin{eqnarray*}
 &a).& \quad H^q(X,\Omega^p_X\otimes\mathcal L)=0,\quad \forall\, p+q>n \\
 &b).& \quad H^q(X,\Omega^p_X\otimes\mathcal L^{-1})=0,\quad \forall\, p+q<n.
\end{eqnarray*}

\subsubsection{Comparison map in the relative case}
The middle dimensional cohomology $
\bigoplus\limits_{p+q=n}H^q(X,\Omega^p_X\otimes \mathcal L^{-1})$
is usually non-zero, and used in the construction of the comparison map connecting deformation Higgs bundle and Hodge theory. Consider a semistable family $f:X\to Y$ over a 1-dimensional base curve $Y$ and denote
$\mathcal L:=\omega_{X/Y}=\Omega^n_{X/Y}(\log \Delta).$
 Given a line bundle $A$ on $Y$ (in the most cases we choose $A$ to be ample) and assume that there is a non-zero section $s$ of $\mathcal L^\nu\otimes f^*A^{-\nu}$ for some $\nu$.
Indeed it is always the case if $f: X\to Y$ is a family of $n$-folds with semi-ample canonical sheaf and with maximal Var$(f)$ and $A$ is a given ample line bundle. By Kawamata-Viehweg's positivity theorem one finds a non-zero section $s$ in $\mathcal L^\nu\otimes f^*A^{-1}$ for $\nu \gg 0$. After replacing the original family by a suitable higher power of the \rev{self-fiber product}
$f^{(r)}: X^{(r)}\to Y$ or by Kawamata base change  $Y'\to Y$ we find a section of $\mathcal L^\nu\otimes f^*A^{-\nu}$ (see \cite{es92}, 3.19 Lemma).

\begin{remark}
	For a family $f:{X}\to{Y}$ of $n$-folds either with good minimal model or of general type.
	Then Kawamata (for good minimal model) and Koll\'ar (for general type)
	showed that $f_*\omega_{{X}/{Y}}^\nu$ is big for $\nu \gg 0$.
	The main difference between the case of good minimal model and the case of semi-ample is that the linear system of $\omega_{{X}/{Y}}^\nu$
	in the first case could be not globally generated over $f^{-1}(U_0)$ for any open subset of $U$,
	while it is globally generated over $f^{-1}(U_0)$ for some open subset of $U$ in the latter case.
	Popa-Schnell \cite{ps17} applied the theory of Hodge module to get a comparison
	similar to what Viehweg-Zuo have done. It has the advantage that one does not care too much about the complication of the singularity appearing in the construction.
	Below, we propose an approach along the original construction by Viehweg-Zuo for a family over a 1-dimensional base curve \cite{vz02}, which works for all above cases and also over higher dimensional bases. We invite the readers to read the details there.
\end{remark}

\begin{proposition}[Viehweg-Zuo]\label{prop-2-3}
 The $\nu$-th cyclic cover defined by a non-zero section $s$ of $\mathcal L^\nu\otimes f^*A^{-\nu}$ induces a family $$g: Z\xrightarrow{\gamma} X\xrightarrow{f}Y$$ with the singular fibers $\Pi$ over $S+T$, where $T$ is the discriminant locus of the "new" singular fibers arising from the cyclic cover $\gamma: Z\to X$. By blowing up of $\Pi$ and we may assume that the reduced singular fibers $\Pi_\text{red}$ is a simple normal crossing divisor.

Taking $(E,\theta)$ to be the graded Higgs bundle of Deligne's quasi-canonical extension of VHS on the middle cohomology $ \mathbb R^ng_*\mathbb Z_{Z\setminus\Pi}$
on $Y\setminus(S+T),$
then there exists a Higgs map $$\rho: (F,\tau)\to (E,\theta)\otimes A^{-1};$$
that is, the following diagram commutes
$$
\xymatrixcolsep{5pc}\xymatrix{
F^{p,q}\ar[r]^-{\tau^{p,q}}\ar[d]^{\rho^{p,q}}& F^{p-1,q+1}\otimes\Omega^1_{{Y}}(\log{S})\ar[d]^{\rho^{p-1,q+1}\otimes\iota}\\
A^{-1}\otimes E^{p,q}\ar[r]^-{\mathrm{id}\otimes\theta^{p,q}}& A^{-1}\otimes E^{p-1,q+1}\otimes\Omega^1_{{Y}}(\log ({S}+{T})).
}
$$
where $\iota: \Omega^1_{{Y}}(\log S)\hookrightarrow \Omega^1_{{Y}}(\log( S+ T))$
 is the natural inclusion.\end{proposition}

 We would like to emphasize the crucial point in the comparison map:
although the Higgs field $\theta$ on $E$ has singularity along $ S+ T$,
its restriction to $\rho(F)$ has only singularity on the original discriminant locus $ S.$

{\bf Sketched proof of \autoref{prop-2-3}.}
Let $D$ denote the zero divisor of $s.$ Note that $D$ could be singular and the intersection of $D$ with the generic fibers could be singular.

{\noindent \bf Step 0. Resolve the singularities. }
By a suitable blowing up
$$\hat f: \hat X\xrightarrow{\sigma}X\to Y,$$
one may assume that $\sigma^*D$
is a normal crossing divisor.
Let $T\subseteq Y$ denote the closure of the discriminant of the map
$$\hat f: \sigma^*D\cap \sigma^{-1}(V)\to U; $$
that is, the locus of $y\in U$ where the simple normal crossing divisor
$\sigma^*D$ meets $\hat f^{-1}(y)$ non-transversally.
 Let
$\Sigma=\hat f^{-1}(T)$, and we take a further blowing up
$$\delta: X'\xrightarrow{\beta} \hat X\xrightarrow{\sigma} X$$
such that $D'+\Delta'+\Sigma':= \delta^{*}(D+\Delta)+\beta^*\Sigma$ is simple normal crossing
and the family
$$f': X'\xrightarrow{\delta}X\xrightarrow{f} Y$$
is log smooth as a morphism between the log pairs
$$ f' :\left(X',(D'+\Delta'+\Sigma')
\right)\to (Y, (S+T)).$$

{\noindent\bf Step 1. Cyclic cover defined by $s$.}
We write $\mathcal M:=\delta^*(\mathcal L\otimes f^*A^{-1})$ and $D':=\delta^*D$, then
$\mathcal M^\nu=O_{X'}(D')$. One takes the $\nu$-th cyclic cover for the section $\delta^*s\in H^0(X',\mathcal M^\nu)$
$$\gamma': Z'\xrightarrow{\text{normalization}}X'(\sqrt[\nu]{\delta^*s})\xrightarrow{\gamma}X'.$$
$Z'$ could be singular.
By taking a resolution of singularity of $ Z'$, and a blowing up at the centers in the fibers
over $Y$ we obtain a non-singular variety $Z$ and a birational map
$\eta: Z\to Z'$. We may assume the induced map
$$g:Z\xrightarrow{\eta}Z'\xrightarrow{\gamma'}X'\xrightarrow{f'}Y$$
is log smooth for the pairs
$$g: (Z , g^{-1}(S+T))\to (Y, (S+T)).$$
We set $\Pi:=g^{-1}(S+T))$, $Z_0=Z\setminus \Pi$ and $\gamma:=\delta\circ \gamma' \circ \eta.$

{\noindent\bf Step 2. Differential forms on the cyclic cover}.
Recall that the local system $\mathbb V=R^ng_*\mathbb Z_0$ over $Y\setminus(S+T)$ gives rise to the filtered logarithmic de Rham bundle
$$\nabla: V\to V\otimes\Omega^1_Y(\log (S+T)),$$
where $\nabla$ is an integrable connection with logarithmic pole along $(S+T),$
as the \rev{quasi-canonical} extension
of $\mathbb V\otimes\mathcal O_{Y\setminus(S+T)}$.
Let $(E,\theta) $ denote the induced system of Hodge bundles
$$\text{Gr}_F(V,\nabla) =(E,\theta)=\Big(\bigoplus_{p+q=n}E^{p,q},\bigoplus_{p+q=n}\theta^{p,q}\Big)$$
with
$$ E^{p,q}=R^qg_*\Omega^p_{Z/Y}(\log\Pi).$$
The Higgs map
$$\theta^{p,q}: E^{p,q}\to E^{p-1,q+1}\otimes\Omega^1_Y(\log (S+T))$$
is
the edge map of $R^\bullet g_*$ of the exact sequence
\begin{equation}\label{eqn-2-3-1}
 0\to g^*\Omega^1_Y(\log(S+T))\otimes\Omega_{Z/Y}^{p-1}(\log\Pi)\to
\Omega^p_Z(\log\Pi)\to
\Omega^p_{Z/Y}(\log\Pi)\to 0.
\end{equation}
We also consider the pulled back of the deformation Higgs bundle $(F,\tau)$ on $Y$ via the blowing up $\delta: X'\to X$
$$\delta^*(F,\tau)=(\bigoplus\delta^*F^{p,q},\bigoplus\delta^*\tau^{p,q})=(\bigoplus F'^{p,q},\bigoplus\tau'^{p,q}),$$
with
$$F'^{p,q}=R^qf'_*(\delta^*\Omega^p_{X/Y}(\log\Delta)\otimes\delta^*\mathcal L^{-1})/\text{torsion}.$$
Note that the Kodaira-Spencer map
$$\tau'^{p,q}: F'^{p,q}\to F'^{p-1,q+1}\otimes\Omega^1_Y(\log S)$$
is the edge map of $R^\bullet f'_*$ of the exact sequence
\begin{equation}\label{eqn-2-3-2}
 0\to f'^*\Omega^1_Y(\log S)\otimes\delta^*\Omega^{p-1}_{X/Y}(\log\Delta)\otimes\mathcal L'^{-1}\to
\delta^*\Omega^p_X(\log\Delta)\otimes\mathcal L'^{-1}\to\delta^*\Omega^p_{X/Y}(\log\Delta)\otimes\mathcal L'^{-1}\to 0.
\end{equation}

{\noindent \bf Step 3. Comparison between deformation Higgs bundle and system of Hodge bundles.}
Let $\bullet$ stand either for $\text{Spec}(\mathbb C)$ or for $Y$. Then the Galois group $\mathbb Z/\nu\mathbb Z$ of
$$\psi: Z\xrightarrow{\eta}Z'\xrightarrow{\gamma'} X'$$ acts on
$\psi_*\Omega^p_{Z/\bullet}(\log\Pi)$
with the eigenspace decomposition
\begin{equation}\label{eqn-2-3-3}
\psi_*\Omega^p_{Z/\bullet}(\log\Pi)
=\Omega^p_{X'/\bullet}(\log(\Delta'+\Sigma'))\oplus
\bigoplus_{i=1}^{\nu-1}\left(\Omega^p_{X'/\bullet}(\log(\Delta'+\Sigma'+D')\otimes\mathcal L'^{-i}\otimes f'^*A^i\right),
\end{equation}
which induces a natural inclusion $\iota$
\begin{equation*}
\xymatrix@R=2cm@C=2.5cm{
\delta^*\Omega_{X/\bullet}^p(\log\Delta)\otimes\mathcal L'^{-1}
\ar@{^(->}[r]^-{\iota}
\ar@{^(->}[d]
&
\psi_*\Omega_{Z/\bullet}^p(\log\Pi)\otimes f'^*(A^{-1})\\
\delta^*\Omega_{X/\bullet}^p(\log\Delta+\Sigma)\otimes\mathcal L'^{-1}
\ar@{^(->}[r]_-{\tiny \txt{Hurwitz formula}}
&
\Omega_{X'/\bullet}^p(\log\Delta'+\Sigma'+D')\otimes\mathcal L'^{-1}\otimes f'^*(A)\otimes f'^*(A^{-1})
\ar@{^(->}[u]^-{}_-{\txt{$1$-th eigen\\ space in (2-6)}}\\
}
\end{equation*}
Via the pulled back $\psi: Z\to X'$ the inclusion $\iota$ together with the natural inclusion
$$\Omega^1_Y(\log S)\hookrightarrow \Omega^1_Y(\log (S+T))$$
induces an inclusion of the exact sequences
$$\psi^* \eqref{eqn-2-3-2}\subseteq \eqref{eqn-2-3-1}\otimes g^*A^{-1},$$ i.e.,
\vspace{1mm}
\begin{equation*} \tiny
\xymatrix@R=5mm@C=3mm{
0
\ar[r]
& \psi^*f'^*\Omega^1_Y(\log S)\otimes\psi^*\delta^*\Omega^{p-1}_{X/Y}(\log\Delta)\otimes\psi^*\mathcal L'^{-1}
\ar[r] \ar@{^(->}[d]
& \psi^*\delta^*\Omega^p_X(\log\Delta)\otimes\psi^*\mathcal L'^{-1}
\ar[r] \ar@{^(->}[d]
& \psi^*\delta^*\Omega^p_{X/Y}(\log\Delta)\otimes\psi^*\mathcal L'^{-1}
\ar[r] \ar@{^(->}[d]
& 0.\\
0
\ar[r]
& g^*\Omega^1_Y(\log(S+T))\otimes\Omega_{Z/Y}^{p-1}(\log\Pi)\otimes g^*A^{-1}
\ar[r]
& \Omega^p_Z(\log\Pi) \otimes g^*A^{-1}
\ar[r]
& \Omega^p_{Z/Y}(\log\Pi)\otimes g^*A^{-1}
\ar[r]
& 0 .
}\vspace{1mm}
\end{equation*}

Finally taking the direct image of the inclusion of the above short exact sequences
$$g_*\left(\psi^* \eqref{eqn-2-3-2}\subseteq \eqref{eqn-2-3-1}\otimes g^*A^{-1}\right),$$
it yields a map between the direct image sheaves
$$\rho^{p,q}: F^{p,q}\to E^{p,q}\otimes A^{-1},$$
which commutes with $\tau$ and $\theta$, as they are just the edge maps connecting the direct image sheaves.
We complete the sketch of the proof of \autoref{prop-2-3}.
\section{ Strict Arakelov Inequality for relative pluri dualizing sheaf of families of manifolds of general type}\label{section-4}
\subsection{ Arakelov Inequality}
In this section we first sketch the proof of the Arakelov inequality \eqref{eqn-1-11} in \autoref{thm-2}.
The main idea in the proof is an application of the general construction performed in \autoref{prop-2-3} to a specific situation and Simpson's Higgs semistability for a system of Hodge bundles.
For reader's convenience we sketch the proof, and the details can be found in \cite{mvz06}.
The proof contains three main steps.

{\noindent \bf $\bullet$ Step I.} Reduce the proof to the case where $W$ is a line subbundle. Indeed, given any non-zero subbundle
$W\subseteq f_*\omega^\nu_{{X}/{Y}}$ of
${\rk} W=r$, by taking the determinant and the $r$-power of \rev{self-fiber product} of $f$,
we have
$$\det W\subseteq (f_*\omega^\nu_{{X}/{Y}})^{\otimes r}\cong \tilde{f}_*\omega^\nu_{\widetilde{X}/{Y}},$$
where $\widetilde{X}$ is the desingularization of the $r$-power of \rev{self-fiber product} ${X} \times_{{Y}} \cdots \times_{{Y}} {X}$
and $\tilde{f}:\, \widetilde{X} \to {Y}$ is the induced fibration.
Hence we may assume $W$ is a line subbundle.

{\noindent\bf $\bullet$ Step II.}
As in \autoref{prop-2-3}, we take the cyclic cover defined by the invertible subsheaf $A\subset f_*\omega^\nu_{X/Y}$ and
 construct a comparison between the deformation Higgs bundle twisted by $A$ and the system of Hodge bundles arising from the cyclic cover.

In order to make this cyclic cover to be possible we replace the original family by suitable base change $Y'\to Y,$ which is
unramified on $U=Y\setminus S$. Such a base change does exist since for $Y=\mathbb P^1$ we assume $\#S\geq 2.$ Hence, we may assume that $A$ is $\nu$-divisible;
that is, there exists an invertible sheaf $A'$ on ${Y}$ such that $A=A'^{\nu}$.
In other words, we get an injection $A'^{\nu} \hookrightarrow f_*\omega^\nu_{{X}/{Y}}$ and hence a non-zero map
$f^*A'^{\nu} \hookrightarrow \omega^\nu_{{X}/{Y}}$.
This is equivalent to a non-zero section $s$ of $\omega^\nu_{{X}/{Y}} \otimes f^*A'^{-\nu}$.
Thus by \autoref{prop-2-3}, we get a new fibration
$$g: (Z,\Pi)\to (Y, S+T), $$
which is log smooth.
Moreover, the graded Higgs bundle $(E,\theta)$ of Deligne's quasi-canonical extension of VHS on the middle cohomology $ R^ng_*\mathbb Z_{Z\setminus\Pi}$
 admits a comparison with the original deformation Higgs bundle $(F,\tau)$ attached to $f:\,X \to Y$;
that is, there exists a Higgs map $$\rho: (F,\tau)\to (E,\theta)\otimes A'^{-1}.$$
It gives the following commutative diagram:
$$
\xymatrixcolsep{5pc}\xymatrix{
F^{p,q}\ar[r]^-{\tau^{p,q}}\ar[d]^{\rho^{p,q}}& F^{p-1,q+1}\otimes\Omega^1_{{Y}}(\log{S})\ar[d]^{\rho^{p-1,q+1}\otimes\iota}\\
A'^{-1}\otimes E^{p,q}\ar[r]^-{\mathrm{id}\otimes\theta^{p,q}}& A'^{-1}\otimes E^{p-1,q+1}\otimes\Omega^1_{{Y}}(\log ({S}+{T})).
}
$$
where $\iota: \Omega^1_{{Y}}(\log S)\hookrightarrow \Omega^1_{{Y}}(\log(S+ T))$
 is the natural inclusion.

{\noindent \bf $\bullet$ Step III.}
The sheaf $A'$ via $\tau$ and $\rho$ generates a Higgs subbundle $$\Big(H=\bigoplus\limits_{q=0}^{n}H^{n-q,q}, \theta|_H\Big) \subseteq (E,\theta),$$
where $H^{n,0}=A'$, and
$$H^{n-q-1,q+1}={\rm Im}\big(\theta|_{H^{n-q,q}}: H^{n-q,q} \to E^{n-q-1,q+1} \otimes \Omega^1_{{Y}}(\log ({S}+{T})) \big) \otimes \Omega^1_{{Y}}(\log{S})^{-1}.$$
Let $q_0\leq n$ be the largest number such that $H^{n-q,q}\neq 0$.
Then
$$\deg H=\sum_{q=0}^{q_0} \deg H^{n-q,q} = \sum_{q=0}^{q_0}\big(\deg A'-q \deg \Omega^1_{{Y}}(\log{S})\big)=(q_0+1)\big(\deg A'-\frac{q_0}{2} \deg \Omega^1_{{Y}}(\log{S})\big).$$
As $(H,\theta)$ is a sub-Higgs bundle of the quasi-canonical extension $(E,\theta)$ of a system of Hodge bundles of a polarized VHS on $Y\setminus (S+T)$ by Simpson's semistability of the quasi-canonical extension of a system of Hodge bundles
 one has $\deg H\leq 0$,
i.e.
$$\deg A =\nu\cdot \deg A'\leq \nu \cdot \frac{q_0}{2} \deg \Omega^1_{{Y}}(\log{S}) .$$
As $q_0\leq n$ and $\deg\Omega^1_Y(\log S)\geq 0,$
we find
$$\deg A\leq \nu \cdot \frac{q_0}{2} \deg \Omega^1_{{Y}}(\log{S})
\leq \frac{n \nu}{2} \deg \Omega^1_{{Y}}(\log{S}).$$
This completes the proof of \autoref{thm-2}.


\subsection{Strict Arakelov Inequality of Families of Varieties of General Type}\label{sec-4-2}
We have seen the Arakelov inequality can be an equality for semistable families of abelian varieties.
In contrast, in this section we shall show \autoref{thm-4} claiming that the Arakelov inequality always holds strictly for families of varieties of general type and for large power $\nu$ such that the relative $\nu$-pluri canonical map is birational.
\begin{proof}[Proof of \autoref{thm-4}]
For simplicity, we prove this for the total direct image sheaf $f_*\omega^\nu_{{X}/{Y}}$;
the proof is similar for subbundles $W\subseteq f_*\omega^\nu_{{X}/{Y}}$ which defines a birational map as in the theorem.

Let $f: X\to Y$ be a family of varieties of general type over a 1-dimensional base $Y$ with discriminant locus $S$
 and assume it is \rev{non-birationally isotrivial}. By \autoref{thm-1} the log curve $(Y,S)$ is hyperbolic. In particular, if $Y=\mathbb P^1$ then $\# S\geq 2.$
Assume on the contrary that there exists such an $\nu \in \mathbb{N}$ satisfying the Arakelov equality
$$\mu (f_*\omega^\nu_{{X}/{Y}})={n\nu\over 2}\cdot\deg\Omega^1_{{Y}}(\log S)=:\mu_0.$$
As for any subbundle $W\subseteq f_*\omega^\nu_{{X}/{Y}}$
by applying \autoref{thm-2} to $W$ we have
$$\mu(W)\leq{n\nu\over 2}\cdot\deg\Omega^1_{{Y}}(\log S)=\mu_0
=\mu (f_*\omega^\nu_{{X}/{Y}} ), $$
i.e., $f_*\omega^\nu_{{X}/{Y}}$ is a semistable vector bundle over ${Y}$.
Consider in the next \rev{step} the $d$-th multiplication map
$$0\to K_{m_d}\to S^d(f_*\omega^\nu_{{X}/{Y}})\xrightarrow{m_d}f_*\omega_{{X}/{Y}}^{d\nu},$$
where $K_{m_d}$ is the kernel of the map $m_d$.
Note that by the definition of the map induced by pluri-canonical linear system the restriction of $K_{m_d}$ on a fiber is the subspace of all \rev{homogeneous} polynomials of degree $d$
in the \rev{homogeneous ideal} defining the birational embedding of that fiber.

Applying again Arakelov inequality in \autoref{thm-2} for the image of the map $m_d$
$$I_{m_d}:= m_d( S^d(f_*\omega^\nu_{{X}/{Y}}) )\subseteq f_*\omega_{{X}/{Y}}^{d\nu}$$ we show
$\mu (I_{m_d})\leq d\cdot\mu_0$.
On the other hand, the symmetric product
$S^d(f_*\omega^\nu_{{X}/{Y}})$ is again semistable of
slope $d\cdot\mu_0$ and $I_{m_d}$ is a quotient bundle we obtain
$$\mu(I_{m_d})\geq \mu\left(S^d(f_*\omega^\nu_{{X}/{Y}})\right)=d\cdot\mu_0,$$
and hence $\mu(I_{m_d})=d\cdot \mu_0.$ From the exact sequence
$$0\to K_{m_d}\to S^d(f_*\omega^\nu_{{X}/{Y}})\to I_{m_d}\to 0$$
we see $\mu(K_{m_d})=d\cdot\mu_0=\mu \left(S^d(f_*\omega^\nu_{{X}/{Y}}) \right)$
and hence $K_{m_d}$ is a semistable subbundle of $S^d(f_*\omega^\nu_{{X}/{Y}})$ of the same slope.

After a base change of ${Y}$ and twisting with a line bundle with a suitable degree we may assume
$f_*\omega^\nu_{{X}/{Y}}$ is semistable of degree zero and
$S^d(f_*\omega^\nu_{{X}/{Y}})$ contains
$K_{m_d}$ as a semistable subbundle of degree zero.
\begin{theorem}[Simpson, cf. \cite{sim92}]
Let $\mathcal C_{dR}$ be the category of vector bundles over ${Y}$ with integrable connections and $\mathcal C_{Dol}$ be the category of semistable Higgs bundle of degree $0$. Then there exists an equivalent functor
$$\mathcal F:\mathcal C_{Dol}\rightarrow\mathcal C_{dR}.$$
\end{theorem}

We just recall some properties about this functor.
Let $(E, 0)$ be a semistable Higgs bundle of degree $0$ with the trivial Higgs field.
Let $(E',0) $ be a sub-Higgs bundle of $(E, 0)$ of degree $0$.
\vspace{2mm}

\noindent (1). The functor $\mathcal F$ preserves the tensor products. In particular it also preserves symmetric powers.

\noindent (2). The underlying bundle of the bundle $\mathcal F((E,0))$ with
the integrable connection is isomorphic to $E$. We call the connection to be canonical and denote it by $\nabla_{can}(E)$.

\noindent (3). The connection $\nabla_{can}(E)$ preserves $E'$ and $\nabla_{can}(E)\mid_{E'}=\nabla_{can}(E')$.

\noindent (4). For a semistable vector bundle $V$ of degree $0$ we may think it is a semistable Higgs bundle with the zero Higgs field. Hence, (1)-(3) above imply that there exists an integrable connection $\nabla$ on $V$ such that for any $d\geq1$
and any subbundle $K\subseteq S^d(V)$ of degree $0$, the connection $S^d(\nabla)$ on $S^d(V)$ preserves $K$.
\vspace{2mm}

Applying (4) for $f_*\omega^\nu_{{X}/{Y}}$ in our situation
we find an integrable connection
$(f_*\omega^\nu_{{X}/{Y}},\nabla) $ such that
$S^d(\nabla)$ preserves $K_{m_d}\subseteq S^d(f_*\omega^\nu_{{X}/{Y}})$ for any $d\in\mathbb N$,
i.e. for each point $p\in U$ we find an analytic open disc
$U_p\subseteq U$ and a flat base
$\mathbb V$ for the solutions of $(f_*\omega^\nu_{{X}/{Y}},\nabla)_{U_p}$
and such that $K_{m_d}\subseteq S^d(f_*\omega^\nu_{{X}/{Y}})$
is spanned by a flat subspace $\mathbb K_{m_d}\subseteq S^d(\mathbb V).$
This means that we find a basis of ${f_*\omega^\nu_{{X}/{Y}}}$ over ${U_p}$ such that under this basis the coefficients of all homogeneous polynomials of degree $d$ in the homogeneous ideal defining the fibers of the family
$f':X':=\eta(X)\to Y$ over $U_p$ are constant,
where $\eta:\,X \dashrightarrow \mathbb{P}_Y^{N}$ with $N=\rank f_*\omega_{X/Y}^{\nu}-1$ is the relative birational embedding defined by $f_*\omega_{X/Y}^{\nu}$.
Hence, we show that the family $f': X'\to Y$ is locally constant over an analytic open discs,
and hence all smooth fibers of $f': X'\to Y$ are isomorphic to each other.
Since $\eta: V\to \eta(V)$ is $U$-birational by the assumption, all fibers of $f:V\to U$ are birational.
By applying \autoref{thm-bbg} due to Bogomolov-B\"ohning-Graf von Bothmer
we show that $f:V\to U$ is birationally isotrivial.
This gives a contradiction since $f$ is \rev{non-birationally isotrivial}.
\end{proof}

\vspace{5mm}
{\noindent \bf Acknowledgment.}
We are grateful to Meng Chen for discussion on the minimal power of pluri-canonical system $|mK_X|$ of three-folds $X$ defining a birational embedding,
and to Yong Hu for a careful reading of an earlier version of our paper and valuable suggestions.
Thanks also to Carlos Simpson for discussion on his \rev{equivalent} functor
between the category of semistable Higgs bundles with trivial Chern classes and the category of vector bundle with integrable connections, and discussion on \autoref{thm-4}.
We also appreciate the anonymous referees for many useful suggestions on
improving the readability of this paper.


\providecommand{\bysame}{\leavevmode\hbox to3em{\hrulefill}\thinspace}
\providecommand{\MR}{\relax\ifhmode\unskip\space\fi MR }
\providecommand{\MRhref}[2]{%
 \href{http://www.ams.org/mathscinet-getitem?mr=#1}{#2}
}
\providecommand{\href}[2]{#2}

\end{document}